\documentclass[11pt,english]{smfart}
\usepackage[T1]{fontenc}

\usepackage[english,french]{babel}
\usepackage{amssymb,url,xspace,smfthm}
\usepackage{commandes}
\usepackage{xcolor}
\usepackage{dsfont}
\usepackage{url}
\usepackage{longtable}

\title[A refinement of Horn's conjecture]{A refinement of Horn's conjecture}
\date{}

\author[A. Médoc]{Antoine Médoc}
\address{IMAG, University of Montpellier, CNRS, Montpellier, France}
\email{antoine.medoc@umontpellier.fr} 

\urladdr{}
\thanks{}

\begin{document}
\begin{abstract}
    We provide a refinement of Horn's conjecture by considering spectra with repetitions. To do this we adapt P. Belkale's techniques to our context, in the form proposed by N. Berline, M. Vergne and M. Walter.
\end{abstract}

\begin{altabstract}
    On propose un affinement de la conjecture de Horn en considérant des spectres avec répétitions. On adapte pour cela les techniques de P. Belkale à notre contexte, dans la forme proposée par N. Berline, M. Vergne et M. Walter.
\end{altabstract}
\maketitle

\tableofcontents

\section{Introduction}

\subsection{Horn's conjecture}

Let $A$ and $B$ be two square matrices of the same order. A natural question (also coming out in physics for example) is to know the relations between the eigenvalues of $A$, $B$ and $A+B$. 
If $A$ and $B$ are diagonalizable and commute, then they are simultaneously diagonalizable and the spectrum of their sum is well known.
In this text we will study the more delicate case of Hermitian matrices with complex coefficients and the related Horn conjecture (proven true in 1999). Interested readers can consult the exposure papers  \cite{fulton,brion,kumar}. A pedagogical introduction can be found in \cite{bhatia}.

Since Hermitian matrices have real eigenvalues, we will see the spectrum of these matrices as tuples with real entries ranked in decreasing order. The previous question can now be reformulated : what are the families $(\Lambda_1, \Lambda_2, \Lambda_3)$ of real tuples such that $\Lambda_1$ (resp. $\Lambda_2$) is the spectrum of an Hermitian matrix $A$ (resp. $B$) and that $\Lambda_3$ is the spectrum of $-(A+B)$ ?
In 1962, Alfred Horn conjectured about the fact that a set of finite inequalities defined by induction are sufficient to describe all possible spectrums for Hermitian matrices and their sum \cite{horn}.

For all $i\in \N^*$ we denote $[i]$ the set of integers $j\in\N$ such that $1\leqslant j\leqslant i$.
Let $r\in\N^*$. We denote by $\R_{\geqslant}^r$ the set of all $\lambda:= (\lambda(i))_{i\in[r]} \in\R^r$ such that $\lambda(1) \geqslant \dots \geqslant \lambda(r)$ and we denote by $\Z^r_\geqslant$ the semi-group $\R^r_\geqslant\cap\Z^r$. 
For all $\lambda\in \R_{\geqslant}^r$ we denote $\OO_{\lambda}$ the set of all hermitian matrices of order $r$ and of spectrum $\lambda$ (this notation comes from the fact that this set is an orbit for the conjugation by the unitary matrices subgroup). 
We will consider an arbitrary number of matrices, not necessarily three : let $s\geqslant 2$ an integer. 
Let's add a variable $t$ to the cartesian product $(\R^r_\geqslant)^s$ by defining
$E(r,s) := (\R_\geqslant^r)^s \x \R.$
We define the Kirwan cone as the set of all $(\Lambda,t) \in E(r,s)$ such that there exists $s$ hermitian matrices with a sum equal to $t\I_r$ ($\I_r$ being the identity matrix of order $r$) and spectrums corresponding to the $s$ real sequences $\Lambda_1,\dots,\Lambda_s$ : 
$$\KKK(r,s) := \ens{ (\Lambda,t) \in E(r,s) \mid t\I_r \in \sum_{l=1}^s \OO_{\Lambda_l} }.$$
We will consider linear equations for the elements of the Kirwan cone $\KKK(r,s)$. 
For all $d\in[r]$ and $J\subset [r]$ subset of cardinality $d$, we identify $J$ with the unique strictly growing map $[d]\rightarrow [r]$ of image $J$ ; we denote, for all $p\in\Z$, $\ga_p(J) := (k-J(k)+p)_{k\in[d]} \in\Z_\geqslant^d$ ; for all $s$-tuple $(\JJ_l)_{l\in[s]}$ of subsets of $[r]$ of cardinality $d$, we denote $\Ga(\JJ) := (\ga_{r-d}(\JJ_l))_{l\in[s]} \in (\Z_\geqslant^r)^s$.

\begin{theo}[Inductive description of the Kirwan cone] \label{kirwanconeinductiontheo}
    Let $(\Lambda,t) \in E(r,s)$. The couple $(\Lambda,t)$ is in $\KKK(r,s)$ if and only if the two following conditions hold :
    \begin{enumerate}
        \item $\sum_{l=1}^s \sum_{j\in [r]} \Lambda_l(j) = rt$ ;
        \item for all $d\in [r-1]$ and all $s$-tuple $(\JJ_l)_{l\in[s]}$ of subsets of $[r]$ of cardinality $d$ such that $(\Ga(\JJ), r-d) \in\KKK(d,s)$, $\sum_{l=1}^s \sum_{j\in \JJ_l} \Lambda_l(j) \leqslant dt$.
    \end{enumerate}
\end{theo}

\begin{proof}
    See subsection \ref{backtohornconjsubsec}.
\end{proof}

This theorem gives less than $2^{rs}+2$ equations to describe the cone $\KKK(r,s)$ and these equations are parameterized by the smallest Kirwan cones. The inductive description of these inequalities is simple enough to allow us to compute them for small dimensions (see section \ref{examplessec}).
The first proof of Horn's conjecture is a consequence of the work of Alexander A. Klyachko \cite{klyachko} and the saturation theorem by Allen Knutson and Terence Tao \cite{kt}. A stronger theorem is presented in \cite{belkale01} and \cite{ktw} (we will see an application in section \ref{examplessec}).

\subsection{Saturation property}

There exists an interesting relation between the Kirwan cone $\KKK(r,s)$ and representation theory \cite{kt}. Let $\U(r)$ be the set of all unitary matrices of order $r$. 
For all $\lambda \in\Z_{\geqslant}^r$ we denote $V(\lambda)$ the irreducible representation of $\U(r)$ with the highest weight $\lambda$. For any representation $V$ of $\U(r)$ we denote $V^{\U(r)}$ the linear subspace of $\U(r)$-invariant vertors of $V$.

\begin{theo}[Knutson-Tao, Saturation property] \label{knutsontaotheorem}
    For all $\Lambda\in (\Z_{\geqslant}^r)^s$, $(\Lambda,0) \in \KKK(r,s)$ if and only if $\parentheses{ \bigotimes_{l=1}^s V(\Lambda_l) }^{\U(r)}\neq\ens0$.
\end{theo}
The original proof of this theorem can be found in \cite{kt} and another can be found in \cite{belkale}. A generalization to quivers was presented in \cite{dw}.

\begin{exem} \label{lrcoeffex}
    We denote $\LR(r,s) := \KKK(r,s)\cap (\Z^r)^s\x\ens 0$.
    The set $\LR(r,3)$ is the semigroup made of all highest weights triplets $(\lambda,\mu, \nu)$ such that the Littlewood-Richardson coefficient $\dim \parentheses{V(\lambda) \otimes V(\mu) \otimes V(\nu)}^{\U(r)}$ is strictly positive.
\end{exem}

\subsection{A refinement of Horn's conjecture} \label{refinementsubsec}

Remark that we are handling many $s$-sequences of elements in the same set. For any set $X$, we consider the natural left action of the symmetric group $\SSS_s$ on the Cartesian product $X^s$ : for all $\s\in\SSS_s$ and $x:=(x_l)_{l\in[s]}\in X^s$ we denote $\s\cdot x := (x_{\s^{-1}(l)})_{l\in [s]}$. Let $\s\in\SSS_s$.

\begin{defi}
    Let $X$ be a set. An element $x:=(x_l)_{l\in[s]}\in X^s$ is \textit{$\s$-stable} if $x = \s\cdot x$. For all subset $A\subset X^s$, the set of all $\s$-stable elements in $A$ is denoted $A^\s$.
\end{defi}

For all $(\Lambda,t)\in E(r,s)$ we denote $\s\cdot(\Lambda, t) := (\s\cdot \Lambda, t)$. Theorem \ref{kirwanconeinductiontheo} gives a set of inequalities to describe the cone $\KKK(r,s)$ and theorem \ref{kirwanconeinductionsymtheo} below assures that a smaller number of these inequalities is enough to describe the cone $\KKK(r,s)^\s$. 
Remark that the second condition in theorem \ref{kirwanconeinductiontheo} using $\Ga$ is well adapted to $\s$-stability : for all $d\in[r]$ and $\JJ:=(\JJ_l)_{l\in[s]}$ $s$-stuple of subsets of $[r]$ of cardinality $d$, $\JJ$ is $\s$-stable if and only if $\Ga(\JJ)$ is $\s$-stable.

\begin{theo} \label{kirwanconeinductionsymtheo}
    Let $(\Lambda,t) \in E(r,s)^\s$. The couple $(\Lambda,t)$ is in $\KKK(r,s)^\s$ if and only if the following conditions hold :
    \begin{enumerate}
        \item $\sum_{l=1}^s \sum_{j\in [r]} \Lambda_l(j) = rt$ ;
        \item for all $d\in [r-1]$ and all $s$-tuple $(\JJ_l)_{l\in[s]}$ of subsets of $[r]$ of cardinality $d$ such that $(\Ga(\JJ), r-d) \in\KKK(d,s)^\s$, $\sum_{l=1}^s \sum_{j\in \JJ_l} \Lambda_l(j) \leqslant dt$.
    \end{enumerate}
\end{theo}

\begin{proof}
    See subsection \ref{consequenceofmaintheosubsec}.
\end{proof}

This is the main result of this paper. Some computations about $\KKK(r,s)$ and $\KKK(r,s)^\s$ are presented in section \ref{examplessec} in the case $s=3$ and $\s = (1~2~3)$ (i.e. with triplets of three equal spectrums).
For example, theorem \ref{kirwanconeinductiontheo} gives $539$ inequalities to describe $\KKK(6,3)$ while theorem \ref{kirwanconeinductionsymtheo} gives only $10$ of them to describe the elements $((\la,\la,\la),t)$ of $\KKK(6,3)$.

\begin{exem} \label{lr63exem}
    Let $\lambda\in\Z^6$. The representation $V(\lambda) \otimes V(\lambda) \otimes V(\la)$ admits a nonzero $U(6)$-invariant vector if and only if $\lambda(1) \geqslant \dots \geqslant \lambda(6)$ and
    $$
\begin{accolade}{rl}
    \lambda(1)+\lambda(2)+\lambda(3)+\la(4)+\la(5)+\la(6) =0,& \\
    \lambda(1)+\lambda(5)+\lambda(6) \leqslant 0,& \\
    \lambda(2)+\lambda(4)+\lambda(6) \leqslant 0,&(*)\\
    \lambda(3)+\lambda(4)+\lambda(5) \leqslant 0.&
\end{accolade}
$$
Remark that inequalitie $(*)$ is in fact a consequence of the others : we can remove $(*)$ from the example. This shows that the number of inequalities given by theorem \ref{kirwanconeinductionsymtheo} is not minimal. We will see how to reduce this number again in Belkale's theorem \ref{belkalethe} and in the main theorem \ref{maintheorem}.
\end{exem}

\subsection*{Acknowledgements}

I would like to thank my doctoral advisor Paul-\'Emile Paradan for his ideas and constant help during the writing of this paper.

\subsection*{Notations and settings}
Most of the notations we use come from \cite{bvw}.
\begin{itemize}
    \item We fix $s\geqslant 2$ (the size of the tuples we study) and $\s\in\SSS_s$ (a permutation preserving the tuples we want to describe).
    \item For all $n,r\in\N^*$ such that $n\geqslant r$, $\Subsets(r,n)$ denote the set of all subsets of $[n]$ made of $r$ elements, which we can identify with the set of all strictly increasing maps $[r] \rightarrow [n]$, and $\Subsets(r,n,s)$ denote the Cartesian product $\Subsets(r,n)^s$.
    \item For all $r,d \in \N^*$ such that $r\geqslant d$, $\JJ \in \Subsets(d,r,s)$ and $\Lambda\in (\R^r)^s$ we denote $T_\JJ(\Lambda) := \sum_{l=1}^s \sum_{j\in\JJ_l} \Lambda_l(j)$ and $T(\Lambda) := \sum_{l=1}^s \sum_{j=1}^r \Lambda_l(j)$ i.e. $T(\JJ) = T_{([r])_{l\in[s]}}(\JJ)$.
    \item In all this paper, $m,d,r,n$ will be positive integers satisfying the three inequalities $m\leqslant d \leqslant r \leqslant n$, $\II$ (resp. $\JJ$) (resp. $\KK$) will be an element of $\Subsets(r,n,s)$ (resp. of $\Subsets(d,r,s)$) (resp. of $\Subsets(m,d,s)$), $U$ will be a complex vector space of finite dimension $n$ and $V$ will be a $r$-dimensional linear subspace of $U$.
    \item We denote by $\Gr(r,U)$ the Grassmannian of all $r$-dimensional linear subspaces of $U$ and we denote by $\Flag(U)$ the set of all complete flags $E:= (E(i))_{i\in[n]}$ on $U$. 
\end{itemize}

\section{Belkale's point of view on theorems \ref{kirwanconeinductiontheo} and \ref{knutsontaotheorem}} \label{bvwviewsec}

\subsection{Belkale's theorem} \label{belkalestheosubsec}

In 2005, Prakash Belkale answered a question from William Fulton \cite{fulton3} and proposed a geometric proof of the Horn conjecture \cite{belkale} using Schubert calculus.
Belkale's geometric point of view is well adapted to prove the refinement presented in this paper.
In 2018, Nicole Berline, Michèle Vergne and Michael Walter presented Belkale's proof \cite{bvw} in a different way and the present text is based on this new redaction.

\subsubsection{Intersecting tuples}

Let $n\in \N^*$, $r\in[n]$ and $I\in \Subsets(r,n)$. Let $U$ be a complex vector space of finite dimension $n$. For all $E\in\Flag(U)$ we denote
$$\Omega_I(E) := \ens{ V\in\Gr(r,U) \mid \forall j\in[r], \rk(V\cap E(I(j))) \geqslant j }$$
the corresponding Schubert subvariety. Its dimension is
$$\dim I := \sum_{j=1}^r (I(j)-j)$$ (see the beginning of section 4 in \cite{fulton} or lemma 3.1.7 in \cite{bvw})
and its class in the integral cohomology ring $H^*(\Gr(r,U))$ is denoted by $\w_I$ ; the dimension $\dim I$ and the cohomology class $\w_I$ only depends on $I$ and does not depend on the flag $E$. 

The class of the point $[\pt]$ is in $H^{2r(n-r)}(\Gr(r,U))$ and, more precisely, $H^{2r(n-r)}(\Gr(r,U)) = \Z [\pt]$. With $a=r(n-r)-\dim I$, $\w_I \in H^{2a}(\Gr(r,U))$.
We will be interested in the product of such cohomology classes. 

\begin{defi} \label{intersectingdef}
    Let $\Intersecting(r,n,s)$ (resp. $\Intersecting^0(r,n,s)$) (resp. $\Intersecting^{00}(r,n,s)$) be the set of all $\II\in \Subsets(r,n,s)$ such that $\prod_{l=1}^s \w_{\II_l}$ is not null (resp. is a multiple of the class of a point) (resp. is the class of a point) in $H^*(\Gr(r,U))$. The elements of $\Intersecting(r,n,s)$ are \textit{intersecting}.
\end{defi}

This is the point of view from \cite{belkale}. In \cite{bvw}, $\II$ is said to be intersecting if, for all $\EE := (\EE_l)_{l\in[s]} \in\Flag(U)^s$, the intersection $$\Omega_\II(\EE) := \bigcap_{l=1}^s \Omega_{\II_l}(\EE_l)$$ is nonempty : by Kleinman's moving lemma, these two definitions are equivalent. The first one can be seen as the most natural one and the second one as the easier one.

\begin{rema} \label{intersectinginclusionsrem}
    Definition \ref{intersectingdef} gives us
    $$\Intersecting^{00}(r,n,s) \subset \Intersecting^0(r,n,s) \subset \Intersecting(r,n,s).$$
\end{rema}

\begin{exem} \label{edimexample}
We have the simple case $$\Subsets(n,n,s) = \Intersecting^{00}(n,n,s)^\s = \ens{([n])_{l\in[s]}}.$$
\end{exem}

\begin{defi}\label{edimdef}
    The \textit{expected dimension} of a tuple $\II\in \Subsets(r,n,s)$ is
    $$\edim \II := r(n-r) - \sum_{l=1}^s (r(n-r) - \dim \II_l) = r(n-r) - T(\Ga(\II)).$$
\end{defi}

A geometrical interpretation of the expected dimension is discussed in subsection 1.1 of \cite{belkale} and in lemma 2.15 of \cite{bvw}. Let $\II\in\Subsets(r,n,s)$, $E\in\Flag(U)$ and $\EE \in \Flag(U)^s$. Since $\Omega_I(E)$ is of dimension $\dim I$, this subvariety is locally described by $\dim \Gr(r,U) - \dim I = r(n-r)-\dim I$ equations. 
Assume that $\II$ is intersecting, i.e. $\Omega_\II(\EE)$ is nonempty.
Let $\CC$ be an irreducible component of $\Omega_\II(\EE)$. It is locally described by $\sum_{l=1}^s (r(n-r) - \dim \II_l)$ equations and its dimension is at least $\edim \II$. 
In fact, there is an equality if the intersection is proper. As we see in definition-proposition \ref{goodsetdefprop}, there is a dense subset $\Good(U,s)$ of $\Flag(U)^s$ such that, if $\EE\in \Good(U,s)$ and $\II$ is intersecting, any irreducible component of $\Omega_\II(\EE)$ is of dimension $\edim\II$.

\begin{rema} \label{intersecting0rem}
    We have $$\Intersecting^0(r,n,s) = \ens{\II \in \Intersecting(r,n,s) \mid \edim \II=0 }.$$
\end{rema}

\subsubsection{The theorem}

Belkale's theorem tells us that a tuple is intersecting (a geometric property) if and only if it satisfies Horn's inequalities (defined by induction). We give a refinement of Belkale's theorem in subsection \ref{refinementofbelkalesubsec}.

For all $d\in[r]$ and $J\in \Subsets(d,r)$, we denote $IJ$ the composition of maps $I\circ J \in\Subsets(d,n)$. 
Let $\II \in \Subsets(r,n,s)$. For all $d\in [r]$ and $\JJ\in\Subsets(d,r,s)$ we denote $\II\JJ := (\II_l \JJ_l)_{l\in [s]}$ which is an element of $\Subsets(d,n,s)$. 

\begin{exem}
    If $I=\ens{2,3,4}$ and $J=\ens{2}$, then $IJ = \ens{3}$.
\end{exem}

\begin{defi} \label{horninequalitiesdefi1}
    Let $\II \in \Subsets(r,n,s)$ and $S\subset \bigcup_{d\in[r]} \Subsets(d,r,s)$. We say that $\II$ verifies the Horn inequalities \eqref{horninequalities1} for the set $S$ if
    \begin{equation} \label{horninequalities1}
        \forall \JJ\in S, \edim \II\JJ  \geqslant \edim\JJ.
    \end{equation}
\end{defi}

\begin{theo}[Belkale] \label{belkalethe}
    Let $\II\in \Subsets(r,n,s)$. The following assertions are equivalent.
    \begin{enumerate}
    \item The tuple $\II$ is intersecting.
    \item We have $\edim \II\geqslant 0$ and, for all $d\in[r-1]$, $\II$ satisfies the Horn inequalities for $\Intersecting(d,r,s)$.
    \item We have $\edim \II\geqslant 0$ and, for all $d\in[r-1]$, $\II$ satisfies the Horn inequalities for $\Intersecting^0(d,r,s)$. 
    \item We have $\edim \II\geqslant 0$ and, for all $d\in[r-1]$, $\II$ satisfies the Horn inequalities for $\Intersecting^{00}(d,r,s)$. 
    \end{enumerate}
\end{theo}

\begin{proof}
    See subsection \ref{belkalesproofsubsec}
\end{proof}

In \cite{belkale}, the fact that the expected dimension is non negative is the codimension condition $(0.1)$. Using example \ref{edimexample}, we could have included this first inequality in the set of the others by taking $d=r$.

Thanks to remark \ref{intersecting0rem}, the third assertion of the theorem allows us to compute $\Intersecting(r,n,s)$ by induction with an easy computation. In \cite{ressayre}, the author presents an inductive algorithm for identifying Littlewood-Richardson coefficients equal to one. In section \ref{examplessec} we will use \cite{buch} to do some computations.

\subsection{Main objects and relations} \label{objectssubsection}

Here we present a very brief resume of the objects used to prove Belkale's theorem \ref{belkalethe}. We are following the presentation made in \cite{bvw}. 
Let $n\in\N^*$ and $U$ a complex vector space of finite dimension $n$.

\subsubsection{Intersecting tuples} Let $r\in[n]$. 

\begin{lemm}\label{intersectingcompositionlemma}
    For all $d\in [r]$, 
    $$\Intersecting(r,n,s) \circ \Intersecting(d,r,s) \subset \Intersecting(d,n,s).$$
\end{lemm}
    This comes from proposition 1 in \cite{fulton2} or from lemma 2.16 in \cite{bvw}.

\begin{lemm} \label{edimposlemma}
    For all $\II\in\Intersecting(r,n,s)$, $\edim \II \geqslant 0.$
\end{lemm}
    This is lemma 4.2.6 in \cite{bvw} which is proven using dominance in algebraic geometry. 

\subsubsection{Harder-Narasimhan lemma on slopes} First, we want to prove that an intersecting tuple satisfies the Horn inequalities for the smaller intersecting tuples.

\begin{defi} \label{schubertposdef}
Let $\EE\in\Flag(U)^s$, $r\in[n]$ and $V\in\Gr(r,U)$. The \textit{Schubert position} of $V$ with respect to $\EE$ is the tuple $\Pos(V,\EE) \in \Subsets(r,n,s)$ such that, for all $l\in[s]$ and $j\in [r]$,
$$\Pos(V,E)_l(j) = \min\ens{ j'\in [r] \mid \dim \EE_l(j')\cap V = j }.$$
For all $\II\in\Subsets(r,n,s)$ we denote 
$$\Omega_\II^0(\EE) := \ens{ V\in \Gr(r,U) \mid \Pos(V,\EE) = \II }.$$ 
\end{defi}

\begin{defi}
Let $r\in[n]$ and $\theta \in (\R^r)^s$. The \textit{slope} associated to $\theta$ is defined by, for all $d\in [r]$ and all $\JJ\in \Subsets(d,r,s)$,
$$\slope_\theta(\JJ) := \frac{1}{d} T_\JJ(\theta) = \frac{1}{d} \sum_{l=1}^s \sum_{j\in \JJ_l} \theta_l(j)\in\R.$$
For all $V\in\Gr(r,U)$, $\FF\in \Flag(V)^s$ and $W\in\Gr(d,V)$ we denote
$\slope_\theta(V,\FF) := \slope_\theta (\Pos(V,\FF)).$
\end{defi}

The Harder-Narasimhan lemma \ref{hnlemma} refers to a classic method used in algebraic geometry. Here, it allows us to compute expected dimensions in a convenient way. Remark that, with the notations $\mu, S, V, \FF, w$ and $\tilde V$ from section 6  of \cite{belkale}, $\mu(S,\tilde V) = -\slope_{-w}(V,\FF)$. 

\begin{lemm}[Harder-Narasimhan] \label{hnlemma}
    Let $r\in[n]$, $V\in\Gr(r,U)$ a linear subspace, $\FF\in\Flag(V)^s$ and $\theta\in (\R^r)^s$ such that, for all $l\in[s]$, $\theta_l$ is increasing. 
    There exists a unique linear subspace $W_*$ of $V$ such that
    $$
    \begin{accolade}{rccl}
        &\slope_\theta(W_*,\FF) &= \min_{W \in \Gr(V), W\neq \ens{0}} \slope_\theta(W,F) &=: m_*,\\
        0< &\dim W_* &= \max_{ W \in \Gr(V), \slope_\theta(W,F) = m_* } \dim W &=: d_*.
    \end{accolade}
    $$
\end{lemm}
This comes from the proof of lemma 4.3.4 from \cite{bvw}.
A version of the Harder-Narasimhan lemma is also used in subsection 6.1 of \cite{belkale}.

\begin{rema} \label{edimslopesrem}
    Let $r\in[n]$ and $d\in[r]$. Let $\II\in \Subsets(r,n,s)$ and $\JJ\in \Subsets(d,r,s)$. 
    Then we have $$\edim \II\JJ - \edim \JJ = d\slope_{-\Ga(\II)} (\JJ) +d(n-r).$$
\end{rema}
It corresponds to the equation following equation $(6.1)$ in \cite{belkale}. This is also lemma 4.3.9 in \cite{bvw}, proven by a direct computation.

\begin{enonce}{Definition-Proposition} \label{goodsetdefprop}
    There is a dense subset $\Good(U,s)$ of $\Flag(U)^s$ which satisfies the following properties.
    \begin{enumerate}
        \item For all $\EE\in\Good(U,s)$ and $r\in[n]$,
        $$\Intersecting(r,n,s) = \ens{ \Pos(V,\EE) ; V\in \Gr(r,U) }.$$
        \item Let $r\in[n]$ and $\II\in \Intersecting(r,n,s)$. For all $\EE,\EE'\in \Good(U,s)$, $\Omega_\II^0(\EE)$ and $\Omega_\II^0(\EE')$ have the same number of irreducible components and each one of them is of dimension $\edim\II$. 
        \item For all $\II\in \Intersecting(r,n,s)$ and $\EE \in \Good(U,s)$, $\Omega_\II^0(\EE)$ is dense in $\Omega_\II(\EE)$.
    \end{enumerate}
\end{enonce}
This comes from propositions 1.1 and 2.3 from \cite{belkale}. It is also proven in lemma 4.3.1 from \cite{bvw}.

Remark that, because of the geometrical interpretation of the expected dimension we discussed after definition \ref{edimdef}, this proposition gives us lemma \ref{edimposlemma} again.

\begin{lemm} \label{horninequalitieslemma}
    Let $r\in[n]$ and $\II\in \Subsets(r,n,s)$ such that, for all $d\in[r]$, $\II$ satisfies the Horn inequalities \eqref{horninequalities1} for $\Intersecting^{00}(d,r,s)$.
    Then, for all $d\in[r]$, $\II$ satisfies the Horn inequalities for $\Intersecting(d,r,s)$.
\end{lemm}
    A proof can be found in section 6.1 from \cite{belkale}. This lemma is also proposition 4.3.10 from \cite{bvw} (with the remark following it). It is a consequence of the Harder-Narasimhan lemma \ref{hnlemma}, remark \ref{edimslopesrem} and definition-proposition \ref{goodsetdefprop} : we will adapt this proof to obtain the refined lemma \ref{horninequalitiessymlemma}.

\begin{lemm}[Horn inequalities] \label{horninequalitieslemma2}
    Let $r\in[n]$, $d\in [r]$ and $\II$ be in $\Intersecting(r,n,s)$. The tuple $\II$ satisfies the Horn inequalities \eqref{horninequalities1} for $\Intersecting(d,r,s)$.
\end{lemm}
This is corollary 4.3.11 in \cite{bvw}. We will adapt this proof to obtain lemma \ref{horninequalitiessymlemma2}.
Remark that lemma \ref{edimposlemma} is a direct consequence of example \ref{edimexample} and lemma \ref{horninequalitieslemma2}. 

\subsubsection{Dimensions and positions of any tuple} \label{tuplesdimensionssubsubsec}

We have already defined the expected dimension of a tuple and we will see two other dimensions : the true dimension in definition \ref{tdimdef} and the kernel dimension in definition \ref{kdimdef}. Both are integers. Finally, we introduce the kernel position of a given tuple in definition \ref{kposdef} : it is another tuple. Let $r\in[n]$ and $\II\in\Subsets(r,n,s)$. 

Let $d\in [r]$. For all $I\in\Subsets(r,n)$ and $J\in\Subsets(d,r)$ we denote 
$$I^J := \ens{ (IJ)(k) - J(k)+k ; k\in[d] } \in\Subsets(d,n-r+d).$$
For all $\JJ\in\Subsets(d,r,s)$ we denote $\II^\JJ := (\II_l^{\JJ_l} )_{l\in[s]}$. 

\begin{lemm} \label{edimformulalemma}
    Let $m\in[d]$, $\JJ\in \Subsets(d,r,s)$ and $\KK\in \Subsets(m,d,s)$. We have
    $$\edim \parentheses{ \II^\JJ \KK } - \edim \KK = \edim \parentheses{ \II\JJ\KK } - \edim (\JJ\KK).$$
\end{lemm}
    This is proven by a direct computation in relation 4.2.10 from \cite{bvw}.

Let $V\in \Gr(r,U)$ and $Q\in \Gr(n-r,U)$ such that $U = V\oplus Q$. We denote
$$\B := \Flag(V)^s\x \Flag(Q)^s.$$
As in definition 2.2 from \cite{belkale}, we denote $\LL(V,Q)$ the $\C$-linear maps from $V$ to $Q$ and, for all $(\FF,\GG)\in \B$, the nonempty set
$$\LL_\II(\FF,\GG) := \bigcap_{l=1}^s \ens{ \ph\in \LL(V,Q) \mid \forall j\in[r], \ph(\FF_l(j)) \subset \GG_l(\II_l(j)-j)}.$$
We also denote 
$$\Pa(\II) := \ens{ (\FF,\GG,\ph) \in \Flag(V)^s\x \Flag(Q)^s\x \LL(V,Q) \mid \ph\in \LL_\II(\FF,\GG) }.$$

\begin{defi} \label{tdimdef}
    The \textit{true dimension} of $\II$ is $\min_{(\FF,\GG)\in \B} \dim \LL_\II(\FF,\GG)$ and denoted
    $\tdim \II$.
\end{defi}

\begin{lemm}\label{tdimedimlemma}
    We have $\tdim\II \geqslant \edim\II$ and the tuple $\II$ is intersecting if and only if this is an equality. 
\end{lemm}
This is the first point of lemma 2.4 and the equivalence $(\alpha) \Leftrightarrow (\ga)$ of proposition 2.3 from \cite{belkale}.

\begin{exem}
    If $\tdim\II=0$ and $\edim \II\geqslant 0$, then $\II$ is intersecting.
\end{exem}

By fixing the true dimension we now consider the nonempty (since $\tdim\II$ is a reached minimum) set 
$$\Pt(\II) := \ens{ (\FF,\GG,\ph) \in \Pa(\II) \mid \dim \LL_\II(\FF,\GG) = \tdim \II }.$$

\begin{defi}\label{kdimdef}
    The \textit{kernel dimension} of $\II$ is 
    $\min_{(F,G,\ph) \in \Pt(\II)} \dim \Ker\ph$ and is denoted $\kdim \II$.
\end{defi}

\begin{lemm} \label{kdimlemma}
    If $\edim \II \geqslant 0$ and $\kdim \II \in\ens{0,r}$, $\II$ is intersecting.
\end{lemm}
This comes from lemma 5.3.4 and corollary 5.2.6 in \cite{bvw}. 

By fixing the kernel dimension we finally consider the nonempty (since $\kdim\II$ is a reached minimum) set 
$$\Pkt(\II) := \ens{ (F,G,\ph)\in\Pt(\II) \mid \dim \Ker \ph = \kdim \II }.$$

\begin{defi} \label{kposdef}
    Assume $\kdim\II \geqslant1$. The \textit{kernel position} of $\II$ is the unique tuple $\kPos(\II) \in \Subsets(\kdim\II, r, s)$ such that, for all $l\in[s]$ and all $k\in [\kdim\II]$,
    $$\kPos(\II)_l(k) = \min_{(\FF,\GG,\ph) \in \Pkt(\II)} \Pos(\Ker\ph, \FF_l)(k).$$
\end{defi}

\begin{lemm} \label{kposintersectinglemma}
    If $1\leqslant \kdim \II \leqslant r-1$, $\kPos(\II)$ is intersecting.
\end{lemm}
This is corollary 5.2.10 in \cite{bvw}.

\begin{rema}
    Using lemma \ref{kposintersectinglemma} and the first point of definition-proposition \ref{goodsetdefprop}, the kernel position of a given tuple is indeed the Schubert position of a linear subspace of $U$ with respect to a sequence of flags on $U$.
\end{rema}

\begin{lemm}[Sherman's relation] \label{recurrencerelationlemma}
    If $1\leqslant \kdim \II \leqslant r-1$,
    $$0\leqslant \tdim \II - \edim \II \leqslant \tdim (\II^{\kPos \II} ) -\edim (\II^{\kPos \II} ).$$
\end{lemm}
    This comes from the last relation of section 3 in \cite{sherman}.
    This is also equation 5.3.3 in \cite{bvw}, proven with topological arguments. 

\begin{rema}
    While the algebraic varieties $\B$, $\Pa(\II)$, $\Pt(\II)$ and $\Pkt(\II)$ depends on $V,Q,U$ (and $\II$), the integers $\edim \II, \tdim\II, \kdim\II$ and the tuple $\kPos(\II)$ only depends on $\II$. 
\end{rema}

\subsection{The proof} \label{belkalesproofsubsec} 
Let (1.), (2.), (3.) and (4.) be the four assertions of Belkale's theorem \ref{belkalethe}.

\subsubsection*{What we have already proven}

Let $r\in\N^*$ and $n\geqslant r$. By lemma \ref{horninequalitieslemma2}, (1.) $\Rightarrow$ (2.). By remarks \ref{intersectinginclusionsrem} and \ref{intersecting0rem}, (2.) $\Rightarrow$ (3.) $\Rightarrow$ (4.). To prove the last implication (4.) $\Rightarrow$ (1.) we follow the proof by induction on $r\in\N^*$ given in \cite{bvw}.

\subsubsection*{Starting the induction}

For all $r\in\N^*$ and $n\geqslant r$, we define $\Horn(r,n,s)$ as the set of all $\II\in \Subsets(r,n,s)$ such that, for all $d\in[r]$, $\II$ satisfies the Horn inequalities \eqref{horninequalities1} for $\Intersecting^{00}(d,r,s)$.
Remark that, because of example \ref{edimexample}, for all $\II \in \Subsets(r,n,s)$, $\edim \II\geqslant 0$ if and only if $\II$ satisfies \eqref{horninequalities1} for $\Intersecting^{00}(r,r,s)$.
For all $r\in \N^*$ we denote by $H(r)$ the assertion
$$\forall d\in [r], \forall n\geqslant d, \Horn(d,n,s) \subset \Intersecting(d,n,s).$$
By example 4.3.12 in \cite{bvw}, $H(1)$ is true. Let $r\geqslant 2$ such that $H(r-1)$ is true. Let $n\geqslant r$. 

\subsubsection*{Contraposition and conclusion}

Let $\II \in \Subsets(r,n,s)$ not intersecting and $d:= \kdim \II$. We want to prove that $\II$ does not satisfy one of the Horn inequalities. If $\edim \II <0$, $\II$ is not in $\Horn(r,n,s)$. Assume that $\edim \II \geqslant 0$.

By lemma \ref{kdimlemma}, $d\in[r-1]$ and we can consider the kernel position $\JJ := \kPos\II$. 
Using lemma \ref{tdimedimlemma} and Sherman's relation \ref{recurrencerelationlemma}, since $\II$ is not intersecting, then $\II^\JJ$ is not intersecting either.
By induction hypothesis $H(r-1)$, since $d\leqslant r-1$, $\II^\JJ$ is not in $\Horn(d,n-r+d,s)$ : there exists $m\in [d]$ and $\KK \in \Intersecting^{00}(m,d,s)$ such that $\edim \II^\JJ \KK < 0$. Using the formula of lemma \ref{edimformulalemma}, 
$$\edim \II(\JJ\KK) - \edim \JJ\KK <0.$$
By lemma \ref{kposintersectinglemma}, $\JJ$ is intersecting. But $\KK$ is also intersecting : by lemma \ref{intersectingcompositionlemma} we deduce that $\JJ\KK$ is intersecting.
Thus, by the contrapositive of lemma \ref{horninequalitieslemma}, there exists $\tilde d\in[r]$ and $\tilde \JJ\in\Intersecting^{00}(\tilde d ,r,s)$ such that $\edim \II \tilde \JJ < 0$ : hence $\II$ is not in $\Horn(r,n,s)$.
From this we deduce $H(r)$ and, by induction, Belkale's theorem \ref{belkalethe}.

\section{A refinement of Belkale's theorem} \label{mainteosection}

\subsection{Main theorem} \label{refinementofbelkalesubsec}

We have introduced the action of the symmetric group $\SSS_s$ in subsection \ref{refinementsubsec}. Let $n\in\N^*$ and $r\in[n]$.
Remark that, for all $\rho \in \SSS_s$, the map
$\papplisansnom{ \II }{ \Subsets(r,n,s)^\s }{ \rho\cdot \II }{ \Subsets(r,n,s)^{\rho \sigma \rho^{-1}} }$ is a bijection.
Let $c_1,\dots,c_p$ be the disjoint cycles of $\s$ ranked by increasing length $l_1\leqslant\dots \leqslant l_p$. We have
$$\Subsets(r,n,s)^\s = \bigcap_{a=1}^p \Subsets(r,n,s)^{c_a}$$
and the only thing interesting about $\s$ is its type $(l_1,\dots, l_p)$ which is a partition of the integer $s = l_1+\dots+l_p$.
From now on, we could replace $\s$ by its conjugate $(1\cdots l_1) \dots (s-l_p\cdots n)$. 
The theorem below is a refinement of Belkale's theorem \ref{belkalethe} which will allow us to prove theorem \ref{kirwanconeinductionsymtheo} using Belkale's method.

\begin{theo} \label{maintheorem}
Let $\II\in \Subsets(r,n,s)^\s$. The following assertions are equivalent.
    \begin{enumerate}
    \item The tuple $\II$ is intersecting.
    \item We have $\edim \II\geqslant 0$ and, for all $d\in[r-1]$, $\II$ satisfies the Horn inequalities \eqref{horninequalities1} for $\Intersecting(d,r,s)^\s$. 
    \item We have $\edim \II\geqslant 0$ and, for all $d\in[r-1]$, $\II$ satisfies the Horn inequalities \eqref{horninequalities1} for $\Intersecting^0(d,r,s)^\s$.
    \item We have $\edim \II\geqslant 0$ and, for all $d\in[r-1]$, $\II$ satisfies the Horn inequalities \eqref{horninequalities1} for $\Intersecting^{00}(d,r,s)^\s$.
    \end{enumerate}
\end{theo}

\begin{proof}
See subsection \ref{maintheoproofsubsec}.
\end{proof}

\begin{exem}
    If $\s= \Id$, theorem \ref{maintheorem} is Belkale's theorem \ref{belkalethe}.
\end{exem}

\begin{exem}
    The set $\Intersecting(5,10,3)$ is made of $718,738$ elements while $\Intersecting(5,10,3)^{(1~2~3)}$ is made of $49$ elements (and none of them is of expected dimension null).
\end{exem}

\subsection{The action of the symmetric group} \label{snactionsubsec}

In this new subsection, we explain how some of the objects and relations of \cite{bvw} (introduced in subsection \ref{objectssubsection}) behave accordingly to the action of $\SSS_s$. 
Let $n\in\N^*$ and $r\in[n]$. Let $\II\in\Subsets(r,n,s)$. 

\subsubsection{In relation with operations on tuples} 

We are interested in $\s$-stability. Let $d\in[r]$ and a tuple $\JJ\in\Subsets(d,r,s)$.

\begin{rema} \label{operationstabilityrem}
    The two operations we have seen on tuples preserves the $\s$-stability : if $\II$ and $\JJ$ are $\s$-stable, $\II\JJ$ and $\II^\JJ$ are $\s$-stable.
\end{rema}

\begin{lemm} \label{symintersectinglemma}
    The tuple $\II$ is intersecting if and only if $\s\cdot \II$ is intersecting.
    We have $\Intersecting(r,n,s)^\s \circ \Intersecting(d,r,s)^\s \subset \Intersecting(d,n,s)^\s$.
\end{lemm}

\begin{proof}
The first assertion comes from the fact that the product is commutative in definition \ref{intersectingdef}.
The second one is a direct consequence of lemma \ref{intersectingcompositionlemma} and remark \ref{operationstabilityrem}.
\end{proof}

\subsubsection{In relation with slopes and the Harder-Narasimhan lemma} We use the uniqueness result of the Harder-Narasimhan lemma \ref{hnlemma} to show that $\s$-stability is preserved.

\begin{rema} \label{slopessymetryrem}
Let $\theta \in (\R^r)^s$. For all $d\in[r]$ and $\JJ\in \Subsets(d,r,s)$, $\slope_\theta(\JJ) = \slope_{\sigma\cdot \theta} (\s\cdot \JJ)$ and $\slope_{-\Ga(\II)}(\JJ) = \slope_{-\Ga(\s\cdot \II)}(\sigma\cdot\JJ).$
\end{rema}

\begin{lemm} \label{hnlemma2}
    Let $\theta$, $\FF$, $d_*$ and $W_*$ as in the Harder-Narasimhan lemma \ref{hnlemma}. Assume that $\FF\in \Good(U,s)$ and that, for all $\JJ\in\Subsets(d_*,r,s)$, we have $\mu_\theta(\JJ) = \mu_\theta(\s\cdot \JJ).$ Then the position $\Pos(W_*,\FF)$ is intersecting and $\s$-stable.
\end{lemm}

\begin{proof}
    Let $\JJ_* := \Pos(W_*,\FF)$. Since $\FF$ is good, $\JJ_*$ is intersecting. In addition to this, $\sigma\cdot \JJ_*$ is also intersecting and there exists $W \in \Gr(d_*,V)$ such that
    $\s\cdot \JJ_* = \Pos(W,F).$
    The linear subspace $W$ is of dimension $d_*$ and, by hypothesis, 
    $$\slope_\theta(W_*,\FF) = \slope_\theta(W,\FF)$$
    By unicity in the Harder-Narasimhan lemma \ref{hnlemma}, $W=W_*$
    hence we have $\JJ_* = \sigma\cdot \JJ_*$.
\end{proof}

\begin{lemm} \label{horninequalitiessymlemma}
    Assume that $\II$ is $\s$-stable and that, for all $d\in[r]$, $\II$ satisfies the Horn inequalities \eqref{horninequalities1} for $\Intersecting^{00}(d,r,s)^\s$.
    Then, for all $d\in[r]$, $\II$ satisfies the Horn inequalities for $\Intersecting(d,r,s)^\s$
\end{lemm}

\begin{proof}
    Using example \ref{edimexample}, the conclusion of the lemma holds for $d=r$.
    Assume that there exists $d\in[r-1]$ and $\JJ \in \Intersecting(d,r,s)^\s$ such that $\edim \II\JJ< \edim\JJ$
    i.e., using remark \ref{edimslopesrem},
    $\slope_{-\Ga(\II)}(\JJ) <-(n-r).$
    
    Let $V\in\Gr(r,U)$ and $\FF\in \Good(V,s)$.
    Since $\JJ$ is intersecting, there exists a nonzero $W\in \Gr(V)$ such that $\JJ = \Pos(W,\FF).$ : 
    $$\slope_{-\Ga(\II)}(W,\FF) < -(n-r).$$
    Using the Harder-Narasimhan lemma \ref{hnlemma}, there exists a unique nonzero $W_*$ of minimal slope $m_*$ with respect to $-\Ga(\II)$ and maximal dimension $d_*$. Remark that, because of the last equation, $m_*<-(n-r).$
    By remark \ref{slopessymetryrem} and lemma \ref{hnlemma2},
    $$\JJ_* := \Pos(W_*,\FF)\in \Intersecting(d_*,r,s)^\s.$$
    For all $W'\in \Omega_{\JJ_*}^0(\FF)$, $\dim W' = d_*$ and $\slope_{-\Ga(\II)}(W',\FF) = m_*$. Since $W_*$ is unique,
    $\Omega_{\JJ_*}^0(\FF) = \ens{W_*}.$
    From this and definition-proposition \ref{goodsetdefprop} we have $\Omega_{\JJ_*}(\FF) = \ens{W_*}$. Since $\FF$ is generic, $\prod_{l=1}^s \w_{(\JJ_*)_l}$ is the class of a point i.e. $\JJ_*\in \Intersecting^{00}(d_*,n,s)$.
    In particular, $\edim \JJ_* = 0$ hence, using remark \ref{edimslopesrem},
    $$\edim \II\JJ_* = d_*m_* + d_*(n-r) <0.$$
    This is in contradiction with the hypothesis on $\II$.
\end{proof}

\begin{lemm}[Horn inequalities] \label{horninequalitiessymlemma2}
    Assume that $\II$ is intersecting and $\s$-stable. Let $d\in [r]$ and $\JJ\in \Intersecting(d,r,s)^\s$. We have $\edim \II \JJ \geqslant \edim \JJ.$
\end{lemm}

\begin{proof}
    For all $d'\in[r]$ and all $\JJ'\in \Intersecting(d',r,s)^\s$, using lemma \ref{intersectingcompositionlemma} $\II\JJ'$ is intersecting and, using lemma \ref{edimposlemma},
    $\edim \II\JJ' \geqslant 0.$
    Particularly, $\II$ satisfies the hypothesis of lemma \ref{horninequalitiessymlemma} and
    $\edim \II \JJ \geqslant \edim \JJ.$
\end{proof}

\subsubsection{In relation with tuples dimensions and the kernel position}

We use the notations of subsubsection \ref{tuplesdimensionssubsubsec}.

\begin{lemm} \label{dimensionsstabilitylemma} Let $\EE\in\Flag(U)^s$, $(\FF,\GG) \in \B$ and $\ph\in\LL(V,Q)$.
\begin{enumerate}
    \item We have $\Omega_{\s\cdot\II}^0(\s\cdot \EE) = \Omega_\II^0(\EE)$ and $\LL_{\s\cdot\II}(\s\cdot \FF,\s \cdot \GG) = \LL_\II(\FF,\GG)$.
    
    \item The triplet $(\FF,\GG,\ph)$ is in $\Pa(\II)$ (resp. $\Pt(\II)$) (resp. $\Pkt(\II)$) if and only if $(\s\cdot \FF,\s\cdot\GG,\ph)$ is in $\Pa(\s\cdot \II)$ (resp. $\Pt(\s\cdot \II)$) (resp. $\Pkt( \s\cdot \II)$).

    \item The three dimensions on tuples we have seen satisfy $\edim (\s\cdot\II) = \edim\II$, $\tdim (\s\cdot \II) = \tdim \II$ and $\kdim (\s\cdot \II) = \kdim \II$.

    \item Assume that $\tdim \II \geqslant 1$. We have $\kPos(\s\cdot \II) = \s\cdot\kPos(\II)$. In particular, if $\II$ is $\s$-stable, then so is $\kPos\II$. 
\end{enumerate}
\end{lemm}

\begin{proof}
We prove the first point of the lemma by reindexing the intersections defining the sets $\Omega_{\s\cdot\II}^0(\s\cdot \EE)$ and $\LL_{\s\cdot\II}(\s\cdot \FF,\s \cdot \GG)$.
In the same way, we prove that the expected dimension is invariant by reindexing the defining sum. In the rest of this proof we consider the now proven equation
\begin{equation}\tag{$*$} \label{listabilityequation}
    \forall (\FF,\GG) \in\B, \LL_{\s\cdot\II}(\s\cdot \FF,\s \cdot \GG) = \LL_\II(\FF,\GG).
\end{equation}

We have $\B = \ens{ (\s^{-1} \cdot \FF, \s^{-1} \cdot \GG) ; (\FF,\GG) \in\B }$. From this and equation \ref{listabilityequation} we know that the true dimensions of $\s\cdot\II$ and $\II$ are the minimum of the same set, hence are equal.

Let $(\FF,\GG) \in\B$ and $\ph\in \LL(V,Q)$. Because of equation \ref{listabilityequation}, $(\FF,\GG,\ph) \in \Pa(\II)$ if and only if $(\s\cdot \FF,\s\cdot \GG,\ph) \in \Pa(\s\cdot\II)$. From this and both equations \ref{listabilityequation} and $\tdim(\s\cdot \II) = \tdim\II$ we deduce that $(\FF,\GG,\ph) \in \Pt(\II)$ if and only if $(\s\cdot \FF,\s\cdot \GG,\ph) \in \Pt(\s\cdot\II)$. Thus the kernel dimensions of $\s\cdot \II$ and $\II$ are the minimum of the same set, hence are equal.

Let $(\FF,\GG) \in\B$ and $\ph\in \LL(V,Q)$. 
We have seen that $\kdim \II = \kdim(\s\cdot\II)$ and that $(\FF,\GG,\ph) \in \Pt(\II)$ if and only if $(\s\cdot \FF,\s\cdot \GG,\ph) \in \Pt(\s\cdot\II)$. Hence $(\FF,\GG,\ph) \in \Pkt(\II)$ if and only if $(\s\cdot \FF,\s\cdot \GG,\ph) \in \Pkt(\s\cdot\II)$. From this we deduce that, for all $l\in[s]$ and $k\in [\kdim\II]$, the integers $\kPos(\II)_{\s^{-1}(l)}(k)$ and $\kPos(\s\cdot\II)_l(k)$ are the minimum of the same set, hence are equal. Finally, $\kPos(\s\cdot \II) = \s\cdot\kPos(\II)$.
\end{proof}

\begin{rema} \label{kpossymrem}
    Using lemma \ref{kposintersectinglemma} and the fourth point of lemma \ref{dimensionsstabilitylemma}, if $\II$ is $\s$-stable and $\kdim\II\geqslant 1$, then $\kPos(\II) \in \Intersecting(\kdim\II,n,s)^\s$. 
\end{rema}

\subsection{Adaptation of Belkale's proof} \label{maintheoproofsubsec}

We prove the main theorem \ref{maintheorem} juste like we proved Belkale's theorem \ref{belkalethe} in subsection \ref{belkalesproofsubsec}, using the point of view of \cite{bvw}.

\subsubsection*{What we have already proven}

Let $r\in\N^*$ and $n\geqslant r$. Lemma \ref{horninequalitiessymlemma2} shows (1.) $\Rightarrow$ (2.). By remarks \ref{intersectinginclusionsrem} and \ref{intersecting0rem}, (2.) $\Rightarrow$ (3.) $\Rightarrow$ (4.). We prove the last implication (4.) $\Rightarrow$ (1.) by induction on $r\in\N^*$.

\subsubsection*{Sarting the induction}

In subsection \ref{snactionsubsec} we have seen the adaptation of the necessary tools to the $\s$-stable case : the main ones are lemma \ref{horninequalitiessymlemma} (playing the role of lemma \ref{horninequalitieslemma}) and remark \ref{kpossymrem} (about the $\s$-stability of the kernel position).
For all $r\in\N^*$ and $n\geqslant r$, we define $\Horn_\s(r,n,s)$ as the set of all $\II\in \Subsets(r,n,s)^\s$ such that, for all $d\in[r]$, $\II$ satisfies the Horn inequalities \eqref{horninequalities1} for $\Intersecting^{00}(d,r,s)^\s$.
Because of example \ref{edimexample}, for all $\II \in \Subsets(r,n,s)$, 
$\edim \II\geqslant 0$ if and only if $\II$ satisfies \eqref{horninequalities1} for $\Intersecting^{00}(r,r,s)^\s$. 
For all $r\in \N^*$ we denote by $H(r)$ the assertion
$$\forall d\in [r], \forall n\geqslant d, \Horn_\s(d,n,s) \subset \Intersecting(d,n,s)^\s.$$
By Belkale's theorem \ref{belkalethe} in the case $r=1$, for all $n\in\N^*$, $\Horn(1,n,s)$ is a subset of $\Intersecting(1,n,s)$ hence, by intersecting with $\Subsets(1,n,s)^\s$, $H(1)$ is true. Let $r\geqslant 2$ such that $H(r-1)$ is true. Let $n\geqslant r$. 

\subsubsection*{Contraposition and conclusion}

Let $\II \in \Subsets(r,n,s)^\s$ not intersecting and $d:= \kdim \II$. We want to prove that $\II$ does not satisfy one of the Horn inequalities.
Just as in subsection \ref{belkalesproofsubsec}, we can suppose that $d \in[r-1]$ and we know that, with $\JJ:= \kPos(\II)$, $\II^\JJ$ is not intersecting.
But, using remarks \ref{kpossymrem} and \ref{operationstabilityrem}, $\II^\JJ$ is $\s$-stable and we can use the induction hypothesis $H(r-1) :$ since $d\leqslant r-1$, $\II^\JJ$ is not in $\Horn_\s(d,n-r+d,s)$. Hence there exists $m\in [d]$ and $\KK \in \Intersecting^{00}(m,d,s)^\s$ such that $\edim \II^\JJ \KK < 0$. Using the formula of lemma \ref{edimformulalemma}, 
$$\edim \II(\JJ\KK) - \edim \JJ\KK <0.$$
By remarks \ref{kpossymrem} and \ref{operationstabilityrem} again, $\JJ\KK$ is $\s$-stable. In addition to this, just like in subsection \ref{belkalesproofsubsec}, $\JJ\KK$ is intersecting.
Thus, by the contrapositive of lemma \ref{horninequalitiessymlemma}, there exists $\tilde d\in[r]$ and $\tilde \JJ\in\Intersecting^{00}(\tilde d ,r,n)^\s$ such that $\edim \II \tilde \JJ < 0$ : hence $\II$ is not in $\Horn_\s(r,n,s)$.
From this we deduce $H(r)$ and, by induction, the main theorem \ref{maintheorem}.

\section{The Kirwan cone with repetitions} \label{kirwanconesection}

\subsection{Back to Horn's conjecture} \label{backtohornconjsubsec}

The link between Belkale's theorem \ref{belkalethe} and Horn's conjecture is presented in subsection 6.3 of \cite{bvw}. We will use it to prove corollary \ref{ktsymcoro}. Let $n\in\N^*$ and $r\in[n]$.

\begin{defi} \label{horninequalitiesdefi2}
    Let $(\Lambda,t)$ and $S\subset \bigcup_{d\in[r]} \Subsets(d,r,s)$. We say that the couple $(\Lambda,t)$ verifies the Horn inequalities \eqref{horninequalities2} for the set $S$ if
    \begin{equation} \label{horninequalities2}
        \forall \JJ\in S, T_\JJ(\Lambda) \leqslant dt.
    \end{equation}
\end{defi}

\begin{theo}[Horn inequalities] \label{horninequalitiesthe}
    For all $(\Lambda,t)\in E(r,s)$, the following assertions are equivalent.
    \begin{enumerate}
        \item The couple $(\Lambda,t)$ is in the Kirwan cone $\KKK(r,s)$.
        \item The tuple $\Lambda$ satisfies $T(\Lambda) = rt$ and, for all $d\in[r-1]$, the Horn inequalities \eqref{horninequalities2} for $\Intersecting(d,r,s)$.
        \item The tuple $\Lambda$ satisfies $T(\Lambda) = rt$ and, for all $d\in[r-1]$, the Horn inequalities \eqref{horninequalities2} for $\Intersecting^{00}(d,r,s)$.
    \end{enumerate}
\end{theo}
This is a consequence of corollary \ref{ktsymcoro}, the proof of which is an adaptation of that of corollary 6.3.3 from \cite{bvw}.

Using theorem \ref{horninequalitiesthe}, if $(\Lambda,t) \in E(r,s)$ satisfies the Horn inequalities \eqref{horninequalities2} for a set of tuples containing $\bigcup_{d=1}^{r-1} \Intersecting^{00}(d,r,s)$, then $(\Lambda,t)$ is in $\KKK(r,s)$. See example \ref{lr3ex} below. 
    In theorem \ref{bktwtheorem} we see that the inequalites parametrized by $\Intersecting^{00}$ are in fact irredudant in the case $s=3$.

\begin{exem} \label{lr3ex}
    Horn's inequalities from theorem \ref{horninequalitiesthe}, Belkale's theorem \ref{belkalethe} and remark \ref{intersecting0rem} give a convenient description of the set $\LR(r,3)$ seen in example \ref{lrcoeffex}. Let $(\lambda,\mu, \nu) \in (\Z_{\geqslant}^s)^3$. The triplet $(\lambda,\mu,\nu)$ is in $\LR(r,3)$ if and only if $T(\lambda,\mu,\nu) =0$ and, for all $d\in[r-1]$ and $\JJ\in \Intersecting^{0}(d,r,3)$, $T_\JJ(\lambda,\mu,\nu) \leqslant 0$. 
\end{exem}

For all $d\in\N$, let $\chi_d := (1)_{k\in[d]}$ be the constant sequence equal to $1$. 

\begin{rema} \label{kirstabilityrem}
    Let $(\Lambda,t) \in E(r,s)$ and $\tau \in\R^s$.
    Remark that $(\Lambda, t) \in \KKK(r,s)$ if and only if
    $$(\Lambda', t') := \parentheses{ \parentheses{ \Lambda_l + \tau_l \chi_r }_{l\in[s]}, t+ \sum_{l\in [s]} \tau_l } \in \KKK(r,s).$$
    Let $S\subset \bigcup_{d\in[r]} \Subsets(d,r,s)$. Remark that $(\Lambda,t)$ satisfies the Horn inequalities \eqref{horninequalities2} for $S$ if and only if $(\Lambda',t')$ satisfies the Horn inequalities \eqref{horninequalities2} for $S$.
\end{rema}

The proof of corollary 6.3.3 in \cite{bvw} also gives us the following result which allows an inductive description on $r\in\N^*$ of $\KKK(r,s)$.

\begin{coro} \label{inductivecor}
    For all $\II \in \Subsets(r,n,s)$, $\II\in \Intersecting^0(r,n,s)^\s$ if and only if $(\Ga(\II),n-r) \in \KKK(r,s)^\s$.
\end{coro}

\begin{proof}
Let $\Ga'(\II) := (\ga_{n-r}(\II_1), \dots, \ga_{n-r}(\II_{s-1}), \ga_0(\II_s)))$. By remark \ref{kirstabilityrem} it is enough to prove that $\II\in \Intersecting^0(r,n,s)$ if and only if $(\Ga'(\II),0) \in \KKK(r,s)$.
Assume that $\II \in\Intersecting^0(r,n,s)$. Using formula 6.3.1 in \cite{bvw} we know that
$$\parentheses{ \bigotimes_{l=1}^s V( \Ga'(\II)_l ) }^{\U(r)} \neq\ens 0.$$ 
Hence, by the Kempf-Ness lemma, $(\Ga'(\II),0) \in\KKK(r,s)$.

We now prove the converse. Assume that $(\Ga'(\II),0) \in\KKK(r,s)$. By Klyachko's corollary 2.13 in \cite{bvw}, it satisfies the Horn inequalities. Using the equations of lemma 4.3.9 in \cite{bvw}, this means that $\II \in \Horn(r,n,s)$. We conclude by Belkale's theorem \ref{belkalethe}.
\end{proof}

Theorem \ref{kirwanconeinductiontheo} is a consequence of Horn's inequalities from theorem \ref{horninequalitiesthe} and of corollary \ref{inductivecor}.

\subsection{A consequence of the refinement of Belkale's theorem} \label{consequenceofmaintheosubsec}
In this subsection we prove theorem \ref{kirwanconeinductionsymtheo} (the inductive description of the $\s$-stable Kirwan cone). Let $n\in\N^*$ and $r\in[n]$.

\begin{coro} \label{ktsymcoro}
    For all $\s$-stable $(\Lambda,t) \in E(r,s)^\s$, the following assertions are equivalent.
    \begin{enumerate}
        \item The couple $(\Lambda,t)$ is in the Kirwan cone $\KKK(r,s)$.
        \item The couple $(\Lambda,t)$ satisfies $T(\Lambda) = rt$ and, for all $d\in[r-1]$, the Horn inequalities \eqref{horninequalities2} for $\Intersecting(d,r,s)$.
        \item The couple $(\Lambda,t)$ satisfies $T(\Lambda) = rt$ and, for all $d\in[r-1]$, the Horn inequalities \eqref{horninequalities2} for $\Intersecting^{00}(d,r,s)^\s$.
    \end{enumerate}
\end{coro}

\begin{proof}
    We adapt the proof of corollary 6.3.3 from \cite{bvw}. 
    Let $(\Lambda,t)$ be a $\s$-stable element of $E(r,s)$ and $\Lambda' := (\Lambda_1, \dots, \Lambda_{s-1}, \Lambda_s - t\chi_r)$. We will use remark \ref{kirstabilityrem}.
    Assume that $(\Lambda,t) \in \KKK(r,s)$ i.e. $(\Lambda',0) \in \KKK(r,s)$ : just as in \cite{bvw} (Klyachko's corollary 2.13), the trace $T(\Lambda')$ is null and, for all $d\in[r-1]$, $(\Lambda',0)$ satisfies the Horn inequalities for all tuples in $\Intersecting(d,r,s)$ hence for all tuples in $\Intersecting(d,r,s)^\s$. We deduce that (1.) $\Rightarrow$ (2.). The second implication (2.) $\Rightarrow$ (3.) is clear.

    To prove that (3.) $\Rightarrow$ (1.), we use the main theorem \ref{maintheorem}.  First we consider the case of integers. 
    Assume that $(\Lambda,t) \in (\Z_{\geqslant}^r)^s\x \Z$, that $T(\Lambda) = rt$ and that, for all $d\in[r-1]$ and $\JJ\in \Intersecting^{00}(d,r,s)^\s$, $T_\JJ(\Lambda) \leqslant dt.$
    Let $p\in\N^*$ be greater than $\max_{l\in [s], j\in [r]} \abs{\Lambda_l(j)} + \abs t$ and 
    $$\tilde \Lambda := \parentheses{ \Lambda_l + p\chi_r }_{l\in [s]} \in (\Z_{\geqslant}^r)^s. $$
    With $n := r+t+sp$ we have, for all $l\in [s]$, $n-r \geqslant \tilde \Lambda_l (1) \geqslant \dots \geqslant \tilde \Lambda_l(r) \geqslant 0$.
    Hence there exists (a unique) $\II \in \Subsets(r,n,s)$ such that $\tilde \Lambda = \Ga(\II)$. 
    As $\Lambda$ is $\s$-stable, so are $\tilde \Lambda$ and $\II$.
    We now want to show that $\II$ is in fact an intersecting tuple of expected dimension $0$. 

    We know that $T(\tilde \Lambda) = T(\Lambda) + rsp = rt+rsp = r(n-r)$ hence, using definition \ref{edimdef}, $\edim \II=0$. For all $d\in[r-1]$ and $\JJ\in \Intersecting^{00}(d,r,s)^\s$, $\edim \JJ =0$ and $T_\JJ(\Lambda) \leqslant dt$ hence, by remark \ref{edimslopesrem},
    $$\edim \II\JJ = -T_\JJ(\tilde \Lambda) +d(n-r) = -T_\JJ(\Lambda) - sdp +d(n-r) \geqslant 0.$$
    Finally $\II$ is $\s$-stable, of expected dimension $0$ and, by the main theorem \ref{maintheorem}, is intersecting. By corollary \ref{inductivecor}, $(\tilde\Lambda, n-r) \in\KKK(r,s)$ and, by remark \ref{kirstabilityrem}, $(\Lambda,t) \in \KKK(r,s)$. 

    We denote by $K_\s(r,s)$ the Klyachko cone made of all $\s$-stable $(\Lambda,t) \in E(r,s)$ such that $T(\Lambda) = rt$ and, for all $d\in[r-1]$, $(\Lambda,t)$ satisfies the Horn inequalities \eqref{horninequalities2} for $\Intersecting^{00}(d,r,s)^\s$.
    At the beginning of the proof we have seen that $\KKK(r,s)^\s\subset K_\s(r,s)$. Then we have proven that $K_\s(r,s) \cap(\Z^r)^s\x \Z$ is a subset of $\KKK(r,s)^\s$. But $\KKK(r,s)$ is invariant under rescaling by any $x\in \R_+$ so $K_\s(r,s) \cap(\Q^r)^s\x \Q \subset \KKK(r,s)^\s$. 
    Since $K_\s(r,s)$ is a polyhedral cone, it is equal to the euclidean closure
    $$\overline{K_\s(r,s) \cap (\Q^r)^s\x \Q} = K_\s(r,s).$$
    On the other hand, $\KKK(r,s)^\s$ is closed for the euclidean topology. Finally, we have the converse inclusion $K_\s(r,s) \subset \KKK(r,s)^\s$ and $\KKK(r,s)^\s = K_\s(r,s)$. 
\end{proof}

\begin{rema} \label{tzsrem}
    As in example \ref{lr3ex}, if $s=3$ and $\s= (1~2~3)$, corollary \ref{ktsymcoro} and the main theorem \ref{maintheorem} gives a convenient description of the set $\LR(r,3)^\s$. Let $\lambda \in \Z_{\geqslant}^r$. The triplet $(\lambda,\lambda,\la)$ is in $\LR(r,3)$ if and only if $\sum_{j\in[r]} \la(j) =0$ and, for all $d\in[r-1]$ and $\JJ := (J,J,J)\in \Intersecting^0(d,r,3)^{\s}$, $\sum_{j\in J} \lambda(j) \leqslant 0$. Finally, together with corollary \ref{inductivecor}, this gives us an inductive description of $\LR(r,3)^\s$. 
\end{rema}

Theorem \ref{kirwanconeinductionsymtheo} is a consequence of corollaries \ref{ktsymcoro} and \ref{inductivecor}.

\section{Examples for $s=3$ and $\s=(1~2~3)$} \label{examplessec}

\subsection{Numbers of equations}

In this subsection we consider the following questions for small values of $r$ : how to describe (with a list of inequalities) the triplets $(\Lambda_1,\Lambda_2,\Lambda_3)$ of real spectra $\Lambda_l(1)\geqslant \dots \geqslant \Lambda_l(r)$ and the real numbers $t$ such that there exists Hermitian matrices $X_1,X_2,X_3$ of order $r$, of spectra $\Lambda_1,\Lambda_2,\Lambda_3$ and of sum $X_1+X_2+X_3=t\I_r$ ? and if $\Lambda_1=\Lambda_2=\Lambda_3$ ?
This correspond to Horn's conjecture with $s=3$ and $\s = (1~2~3)$. In the tabular below, we compute some of the following numbers using the main theorem \ref{maintheorem} and \cite{buch}.
\begin{itemize}
    \item The integer $l^{0}(r,3)$ is the number of inequalities given by theorem \ref{kirwanconeinductiontheo} to describe $\KKK(r,3)$ : the $3(r-1)$ inequalities $\Lambda_l(i)\geqslant \Lambda_l(i+1)$ (for all $l\in [3]$ and $i\in [r-1]$) ; the two inequalities $T(\Lambda) \leqslant rt$ and $T(\Lambda) \geqslant rt$ ; the inequalities $T_{\JJ} \leqslant dt$ for all $d\in[r-1]$ and $\JJ\in \Intersecting^0(d,r,3)$.
    \item The integer $l_{\min}(r,3)$ is the minimal number of inequalities taken from the previous ones to describe the cone $\KKK(r,3)$. Using the Belkale-Knutson-Tao-Woodward theorem \ref{bktwtheorem}, if $r\geqslant 2$ we remove inequalities $T_\JJ(\Lambda) \leqslant dt$ for all $d\in[r-1]$ and $\JJ \notin \Intersecting^{00}(d,r,3)$.
    If $r=2$, we also remove the three inequalities $\Lambda_l(1)\geqslant \Lambda_l(2)$ (for all $l\in [3]$).
    \item The integer $l_\s^0(r,3)$ is the number of inequalities given by theorem \ref{kirwanconeinductionsymtheo} to describe the cone $\KKK(r,3)^\s$, i.e. the element of the form $((\lambda,\lambda,\lambda),t)$ in $\KKK(r,3)$ : the $r-1$ inequalities $\lambda(i)\geqslant \lambda(i+1)$ (for all $i\in[r-1]$) ; the two inequalities $3\sum_{j\in[r]} \lambda(j) \leqslant rt$ and $3\sum_{j\in[r]} \lambda(j) \geqslant rt$ ; the inequalities $3 \sum_{j\in J} \lambda(j) \leqslant dt$ for all $d\in[r-1]$ and $(J,J,J) \in\Intersecting^0(d,r,3)$. Clearly, $l_\s(r,3) \leqslant l(r,3) = l_{\Id}(r,3)$.
    \item The integer $l_\s^{00}(r,3)$ is the smaller number of inequalities we found in corollary \ref{ktsymcoro} to describe the cone $\KKK(r,3)^\s$. We remove from the previous inequalities the ones of the form $3 \sum_{j\in J} \lambda(j) \leqslant dt$ for all $d\in[r-1]$ and $(J,J,J) \notin \Intersecting^{00}(d,r,3)$. We do not know if they are irredudant.
\end{itemize}

\begin{center}
\begin{tabular}{|c|c|c|c|c|c|c|c|c|c|c|}
    \hline
     $r$&$1$&$2$&$3$&$4$&$5$&$6$&$7$&$8$&$9$&$10$  \\ \hline
     $l^0(r,3)$       & $2$ & $8$ & $20$ & $52$ & $156$ & $539$ & $2,082$ & $8,775$ & $39,742$ & $191,382$ \\ \hline
     $l_{\min}(r,3)$  & $2$ & $5$ & $20$ & $52$ & $156$ & $538$ & $2,062$ & $8,522$ & $37,180$ & $168,602$ \\ \hline
     $l_\s^0(r,3)$      & $2$ & $3$ & $4$  & $7$  & $10$  & $10$  & $18$    & $25$    & $24$     & $51$      \\ \hline
     $l^{00}_\s(r,3)$ & $2$ & $3$ & $4$  & $7$  & $10$  & $9$   & $16$    & $21$    & $18$     & $35$      \\ \hline
\end{tabular}
\end{center}

\begin{exem} \label{0=00exem}
    For all $d,r\in\N$ such that $1\leqslant d \leqslant r \leqslant 5$, $$\Intersecting^0(d,r,3) = \Intersecting^{00}(d,r,3).$$
\end{exem}

\begin{theo}[Belkale-Knutson-Tao-Woodward] \label{bktwtheorem}
    Let $(\Lambda,t) \in E(r,s)$. The couple $(\Lambda,t)$ is in $\KKK(r,3)$ if and only if the following conditions hold :
    \begin{enumerate}
        \item for all $l\in[3]$ and $i\in [r-1]$, $\Lambda_l(i) \geqslant \Lambda_l(i+1)$ ;
        \item $T(\Lambda) = rt$ (seen as two inequalities) ;
        \item if $r\geqslant 3$, for all $d\in [r-1]$ and all $\JJ\in \Intersecting^{00}(d,r,3)$, $T_\JJ(\Lambda) \leqslant dt$.
    \end{enumerate}
    In addition to this, all of these inequalities are essential. 
\end{theo}

\begin{proof}
    The fact that these inequalities are enough to describe $\KKK(r,s)$ comes from theorem \ref{horninequalitiesthe}. The fact that these inequalities are essential is theorem 4 from \cite{ktw}.
\end{proof}

\subsection{Equations}

For all $r,d\in\N$ such that $1\leqslant d<r\leqslant 5$ we give the list of the elements $(\JJ_1,\JJ_2,\JJ_3) \in \Intersecting^{00}(d,r,3)$ (remember example \ref{0=00exem}) up to permutation ; the $\s$-stable ones are in bold. Other examples can be found in appendix A from \cite{bvw} and in subsection 1.3 from \cite{klyachko} who shows how some historical equations are consequences of theorem \ref{horninequalitiesthe}.

\begin{exem}
    Description of $\Intersecting^{00}(d,2,3)$ for $d\in[1]$.
    \begin{center}
    \begin{tabular}{|c|ccc|}
    \hline
    $d$ & $\JJ_1$ & $\JJ_2$ & $\JJ_3$ \\ \hline
    $1$& $\ens1$ & $\ens2$ & $\ens2$ \\
    \hline
    \end{tabular}
    \end{center}
\end{exem}

\begin{exem}
    Description of $\Intersecting^{00}(d,3,3)$ for $d\in[2]$.
    \begin{center}
    \begin{tabular}{|c|ccc|}
    \hline
    $d$ & $\JJ_1$ & $\JJ_2$ & $\JJ_3$ \\ \hline
    $1$& $\ens1$ & $\ens3$ & $\ens3$ \\
    & $\ens 2$ & $\ens 2$ & $\ens3$\\ \hline
    $2$ & $\ens{1,2}$ & $\ens{2,3}$ & $\ens{2,3}$\\
    & $\ens{1,3}$ & $\ens{1,3}$ & $\ens{2,3}$\\
    \hline
    \end{tabular}
    \end{center}
\end{exem}

\begin{exem}
    Description of $\Intersecting^{00}(d,4,3)$ for $d\in[3]$.
    \begin{center}
    \begin{longtable}{|c|ccc|}
    \hline
    $d$ & $\JJ_1$ & $\JJ_2$ & $\JJ_3$ \\ \hline
    $1$& $\ens1$ & $\ens4$ & $\ens4$ \\
    & $\ens2$ & $\ens3$ & $\ens4$ \\
    & $\mathbf{\ens3}$ & $\mathbf{\ens3}$ & $\mathbf{\ens3}$ \\\hline
    $2$ & $\ens{1,2}$ & $\ens{3,4}$ & $\ens{3,4}$\\
     & $\ens{1,3}$ & $\ens{2,4}$ & $\ens{3,4}$\\
     & $\ens{1,4}$ & $\ens{1,4}$ & $\ens{3,4}$\\
     & $\ens{1,4}$ & $\ens{2,4}$ & $\ens{2,4}$\\
     & $\ens{2,3}$ & $\ens{2,3}$ & $\ens{3,4}$\\
     & $\ens{2,3}$ & $\ens{2,4}$ & $\ens{2,4}$\\\hline
     $3$ & $\ens{1,2,3}$ & $\ens{2,3,4}$ & $\ens{2,3,4}$\\
      & $\ens{1,2,4}$ & $\ens{1,3,4}$ & $\ens{2,3,4}$\\
      & $\mathbf{\ens{1,3,4}}$ & $\mathbf{\ens{1,3,4}}$ & $\mathbf{\ens{1,3,4}}$\\
    \hline
    \end{longtable}
    \end{center}
\end{exem}

\begin{exem}
    Description of $\Intersecting^{00}(d,5,3)$ for $d\in[4]$.
    \begin{center}
    \begin{longtable}{|c|ccc|}
    \hline
    $d$ & $\JJ_1$ & $\JJ_2$ & $\JJ_3$ \\ \hline
    $1$& $\ens1$ & $\ens5$ & $\ens5$ \\
    & $\ens2$ & $\ens4$ & $\ens5$ \\
    & $\ens3$ & $\ens3$ & $\ens5$ \\
    & $\ens3$ & $\ens4$ & $\ens4$ \\\hline
    $2$ & $\ens{1,2}$ & $\ens{4,5}$ & $\ens{4,5}$\\
     & $\ens{1,3}$ & $\ens{3,5}$ & $\ens{4,5}$\\
     & $\ens{1,4}$ & $\ens{2,5}$ & $\ens{4,5}$\\
     & $\ens{1,4}$ & $\ens{3,5}$ & $\ens{3,5}$\\
     & $\ens{1,5}$ & $\ens{1,5}$ & $\ens{4,5}$\\
     & $\ens{1,5}$ & $\ens{2,5}$ & $\ens{3,5}$\\
     & $\ens{2,3}$ & $\ens{3,4}$ & $\ens{4,5}$\\
     & $\ens{2,3}$ & $\ens{3,5}$ & $\ens{3,5}$\\
     & $\ens{2,4}$ & $\ens{2,4}$ & $\ens{4,5}$\\
     & $\ens{2,4}$ & $\ens{2,5}$ & $\ens{3,5}$\\
     & $\ens{2,4}$ & $\ens{3,4}$ & $\ens{3,5}$\\
     & $\mathbf{\ens{2,5}}$ & $\mathbf{\ens{2,5}}$ & $\mathbf{\ens{2,5}}$\\
     & $\ens{2,5}$ & $\ens{2,5}$ & $\ens{3,4}$\\
     & $\mathbf{\ens{3,4}}$ & $\mathbf{\ens{3,4}}$ & $\mathbf{\ens{3,4}}$\\\hline
     $3$ & $\ens{1,2,3}$ & $\ens{3,4,5}$ & $\ens{3,4,5}$\\
         & $\ens{1,2,4}$ & $\ens{2,4,5}$ & $\ens{3,4,5}$\\
         & $\ens{1,2,5}$ & $\ens{1,4,5}$ & $\ens{3,4,5}$\\
         & $\ens{1,2,5}$ & $\ens{2,4,5}$ & $\ens{2,4,5}$\\
         & $\ens{1,3,4}$ & $\ens{2,3,5}$ & $\ens{3,4,5}$\\
         & $\ens{1,3,4}$ & $\ens{2,4,5}$ & $\ens{2,4,5}$\\
         & $\ens{1,3,5}$ & $\ens{1,3,5}$ & $\ens{3,4,5}$\\
         & $\ens{1,3,5}$ & $\ens{1,4,5}$ & $\ens{2,4,5}$\\
         & $\ens{1,3,5}$ & $\ens{2,3,5}$ & $\ens{2,4,5}$\\
         & $\mathbf{\ens{1,4,5}}$ & $\mathbf{\ens{1,4,5}}$ & $\mathbf{\ens{1,4,5}}$\\
         & $\ens{1,4,5}$ & $\ens{2,3,5}$ & $\ens{2,3,5}$\\
         & $\ens{2,3,4}$ & $\ens{2,3,4}$ & $\ens{3,4,5}$\\
         & $\ens{2,3,4}$ & $\ens{2,3,5}$ & $\ens{2,4,5}$\\
         & $\mathbf{\ens{2,3,5}}$ & $\mathbf{\ens{2,3,5}}$ & $\mathbf{\ens{2,3,5}}$\\ \hline  
    $4$ & $\ens{1,2,3,4}$ & $\ens{2,3,4,5}$ & $\ens{2,3,4,5}$ \\
        & $\ens{1,2,3,5}$ & $\ens{1,3,4,5}$ & $\ens{2,3,4,5}$ \\
        & $\ens{1,2,4,5}$ & $\ens{1,2,4,5}$ & $\ens{2,3,4,5}$ \\
        & $\ens{1,2,4,5}$ & $\ens{1,3,4,5}$ & $\ens{1,3,4,5}$ \\
    \hline
    \end{longtable}
    \end{center}
    Let $\lambda \in \R^5$. Using corollary \ref{ktsymcoro} and this tabular we know that $(\lambda, \lambda, \lambda)$ is in $\LR(5,3)$ (see example \ref{lrcoeffex}) if and only if $\lambda(1) \geqslant \dots \geqslant \la(5)$ and
    $$\begin{accolade}{rl}
    \lambda(1)+\lambda(2)+\lambda(3)+\la(4)+\la(5)& =0\\
    \lambda(2) + \lambda(5) &\leqslant 0\\
    \lambda(3)+\lambda(4) &\leqslant 0\\
    \lambda(1)+\lambda(4)+\lambda(5) &\leqslant 0\\
    \lambda(2)+\lambda(3)+\lambda(5) &\leqslant 0
    \end{accolade}.$$
\end{exem}

\begin{exem}
    In example \ref{lr63exem}, we can also see that equation $(*)$ is the consequence of the others because $(\ens{2,4,6}, \ens{2,4,6}, \ens{2,4,6})$ is in $\Intersecting^0(3,6,2)^\s$ but not in $\Intersecting^{00}(3,6,2)^\s$
\end{exem}

\bibliography{bibtemplate}
\bibliographystyle{smfalpha}

\end{document}